\documentclass{amsart}
\usepackage{amsmath,amsfonts,amssymb,amscd,amsthm,amsbsy,epsf}
\newtheorem{thm}{Theorem} 
\newtheorem{cor}{Corollary} 

\newtheorem*{Pittsinequality}{Pitt's inequality}
\newtheorem*{Cor-Hardy}{Corollary 1}
\newtheorem*{Cor-Young}{Corollary 2}
\newtheorem{lem}{Lemma}
\newtheorem*{SW}{Stein-Weiss Lemma}
\newtheorem*{SWthm}{Stein-Weiss Theorem}
\theoremstyle{definition}
\newtheorem*{remark}{Remark}
\def\F{{\mathcal F}}
\def\S{{\mathcal S}}
\def\real{{\mathbb  R}}

\def\step#1{\medskip\noindent$\underline{\text{Step #1}}$.\enspace}

\def\smallnegint{\mathop{\int\mkern-13mu
        \raise.5ex\hbox{${\scriptscriptstyle\diagup}$}}\nolimits}

\begin{document}

\title{Pitt's inequality with sharp convolution estimates}
\author{William Beckner}
\address{Department of Mathematics, The University of Texas at Austin,
1 University Station C1200, Austin TX 78712-0257 USA}
\email{beckner@math.utexas.edu}
\begin{abstract}
Sharp $L^p$ extensions of Pitt's inequality expressed as a weighted
Sobolev inequality are obtained using convolution estimates and Stein-Weiss
potentials. Optimal constants are obtained for the full Stein-Weiss 
potential as a map from $L^p$ to itself which in turn yield semi-classical
Rellich inequalities on $\real^n$. Additional results are obtained for
Stein-Weiss potentials with gradient estimates and with mixed homogeneity.
New proofs are given for the classical Pitt and Stein-Weiss estimates.
\end{abstract}
\maketitle

Weighted inequalities for the Fourier transform provide a natural measure
to characterize both uncertainty and the balance between functional growth
and smoothness. On $\real^n$ the question is to determine sharp quantitative
comparisons between the relative size of a function and its Fourier 
transform
at infinity. Pitt's inequality illustrates this principle at the spectral
level (see [4], [7]):
\begin{equation}\label{eq1}
\int_{\real^n} \Phi (1/|x|) |f(x)|^2\,dx
\le C_\Phi \int_{\real^n} \Phi (|y|) |\hat f(y)|^2\,dy
\end{equation}
where $\Phi$ is an increasing function, $f$ is in the Schwartz class
$\S(\real^n)$ and the Fourier transform is defined by
\begin{equation*}
(\F f) (y) = \hat f(y) = \int_{\real^n} e^{2\pi ixy} f(x)\,dx\ .
\end{equation*}
The objective here is to extend Pitt's inequality to $L^p (\real^n)$ for
$1<p<\infty$ in the form:
\begin{equation}\label{eq2}
\int_{\real^n} |\Phi (1/|x|) f(x) |^p \,dx
\le C_{\Phi,p} \int_{\real^n} |\Phi (\sqrt{-\Delta}\,)f(x)|^p \,dx
\end{equation}
where $\Delta$ denotes the Laplacian on $\real^n$.
Sharp constants for Pitt's inequality can be calculated using Stein-Weiss
potentials and Young's inequality as in the equivalent results ([4]):

\subsection*{Pitt's inequality}
For $f\in \S(\real^n)$ and $0\le \alpha <n$
\begin{gather}
\int_{\real^n} |x|^{-\alpha}|f(x)|^2 \,dx
\le C_\alpha \int_{\real^n} |y|^\alpha |\hat f(y)|^2\,dy
\label{eq3}\\
\noalign{\vskip6pt}
C_\alpha = \pi^\alpha \left[ \Gamma \Big(\frac{n-\alpha}4\Big) \Big/
\Gamma \Big(\frac{n+\alpha}4\Big) \right]^2 \notag
\end{gather}

\subsection*{Stein-Weiss potentials}
For $f\in L^2 (\real^n)$ and $0<\alpha< n$
\begin{gather}
\Big| \int_{\real^n\times\real^n} f(x) \frac1{|x|^{\alpha/2}}\
\frac1{|x-y|^{n-\alpha}}\ \frac1{|y|^{\alpha/2}} f(y) \,dx\,dy \Big|
\le B_\alpha \int_{\real^n}  |f(x)|^2\,dx 
\label{eq4}\\
\noalign{\vskip6pt}
B_\alpha = \pi^{n/2}
\left[\Gamma\Big(\frac{\alpha}2\Big) \Big/
\Gamma \Big(\frac{n-\alpha}2\Big)\right]
\left[\Gamma\Big(\frac{n-\alpha}4\Big) \Big/
\Gamma \Big(\frac{n+\alpha}4\Big)\right]^2\ .\notag
\end{gather}
Alternate proofs of inequality (3) are given in [11] and [24].
This Stein-Weiss potential corresponds to a linear operator in $L^2 
(\real^n)$
\begin{equation*}
g\to |x|^{-\alpha/2}\ (|x|^{-(n\,-\,\alpha/2)} * g)
\end{equation*}
and inequality \eqref{eq4} can be rephrased
\begin{equation}\label{eq5}
\Big\|\, |x^{-\alpha/2} \big(|x|^{-(n\,-\,\alpha/2)} *g\big)
\Big\|_{L^2 (\real^n)}
\le \pi^{n/2}
\left[ \frac{\Gamma (\frac{\alpha}4)\ \Gamma (\frac{n-\alpha}4)}
{\Gamma (\frac{2n-\alpha}4)\ \Gamma (\frac{n+\alpha}4)} \right]
\|g\|_{L^2(\real^n)}
\end{equation}
where $0<\alpha <n$.
It is this inequality that extends  naturally to $L^p(\real^n)$ using
the Hardy-Littlewood paradigm that a positive integral operator that 
commutes
with dilations may be reduced to Young's inequality for sharp convolution
estimates on the multiplicative group
(see discussion on bilinear forms and integrals in \cite{14}, chapter 9).
This method has been applied by the author in several articles, including
\cite4, \cite5 and \cite7, but it was used earlier by Herbst in his study
of the Klein-Gordon equation for a Coulomb potential (see Theorem 2.5 and
the related discussion in \cite{15}).

\begin{thm}\label{thm1}
{\bf (Herbst)}
For $g\in L^p(\real^n)$ and $h\in L^{p'} (\real^n)$ 
with $1<p<\infty$, $1/p + 1/p' =1$ and $0<\alpha <n$
\begin{gather}
\Big\| \, |x|^{-\alpha/p} \big(|x|^{-(n\,-\,\alpha/p)} *g\big)
\Big\|_{L^p(\real^n)}  
\le C_{\alpha,p} \|g\|_{L^p(\real^n)} \label{eq6}\\
\noalign{\vskip6pt}
\Big\|\, |x|^{-(n\,-\,\alpha/p)} * \big( |x|^{-\alpha/p} h\big)
\Big\|_{L^{p'} (\real^n)}
\le C_{\alpha,p} \|h\|_{L^{p'} (\real^n)} \label{eq7}\\
\noalign{\vskip6pt}
C_{\alpha,p} = \pi^{n/2} \left[ \frac{\Gamma (\frac{\alpha}{2p})
\Gamma (\frac{n-\alpha}{2p}) \Gamma (\frac{n}{2p'})}
{\Gamma (\frac{n}2 - \frac{\alpha}{2p})
\Gamma (\frac{n}{2p'} + \frac{\alpha}{2p})
\Gamma (\frac{n}{2p})} \right]\ .
 \notag
\end{gather}
\end{thm}

\begin{proof}
Observe that \eqref{eq7} follows from \eqref{eq6} by duality.
By using rearrangement and symmetrization, both inequalities are reduced
to non-negative radial decreasing functions.
Set $t=|x|$, $k(t) = |x|^{n/p} g(|x|)$, $\sigma = \frac{n}p - \frac{n}2
-\frac{\alpha}{2p}$ and
\begin{equation*}
\psi_\sigma (t) =   t^\sigma
\int_{S^{n-1}} \Big( t + \frac1t - 2\xi_1\Big)^{-(\sigma +n/p')} \,d\omega
\end{equation*}
where $d\omega$ is  surface measure on $S^{n-1}$ and $\xi_1$ is
the first component of $\xi \in S^{n-1}$.
Note that $\sigma$ can take both positive and negative values, 
but $\sigma + \frac{n}{p'} = \frac{n}2 -\frac{\alpha}{2p} >0$ and 
$\psi_\sigma \in L^r (\real_+)$ for $r\ge1$.
Then inequality \eqref{eq6} corresponds to
\begin{equation}\label{eq8}
\|\psi_\sigma * k \|_{L^p(\real_+)}
\le \|\psi_\sigma\|_{L^1(\real_+)} \|k\|_{L^p(\real_+)}
\end{equation}
with $C_{\alpha,p} = \|\psi_\sigma \|_{L^1(\real_+)}$.
\begin{equation*}
\begin{split}
\|\psi_\sigma\|_{L^1(\real_+)}
& = \int_0^\infty \psi_\sigma (t) \frac1t \,dt
= \int_{\real^n} |x-\xi|^{-(n\,-\, \alpha/p)} 
|x|^{-(n/p'\,+\,\alpha/p)}\,dx\\
\noalign{\vskip6pt}
& = \pi^{n/2} \left[ \frac{\Gamma (\frac{\alpha}{2p} )
\Gamma (\frac{n-\alpha}{2p}) \Gamma (\frac{n}{2p'})}
{\Gamma (\frac{n}2 - \frac{\alpha}{2p})
\Gamma (\frac{n}{2p'}+ \frac{\alpha}{2p})
\Gamma (\frac{n}{2p})} \right]
\end{split}
\end{equation*}
where $|\xi| =1$.
\renewcommand{\qed}{}
\end{proof}

As outlined in the introduction, this result can be expressed as a
weighted Sobolev inequality.
Because such estimates have the same dilation character as Hardy's
inequality, they have been classified as Hardy-Rellich inequalities and
in addition they can be viewed in the framework of the
Maz'ya-Eilertsen inequality \cite{11}.
The following corollary to Theorem~\ref{thm1} generalizes
equation~(20) in \cite7.
The techniques developed here allow the extension of the results
from \cite7 on iterated Stein-Weiss potentials from
$L^2(\real^n)$ to $L^p(\real^n)$ (see also Corollary 14 in \cite{10}).

\begin{cor}\label{cor1} 
For $g\in \S(\real^n)$, $1<p<\infty$, $0\le \alpha <n$ 
\begin{equation}\label{eq9}
\big\|\, |x|^{-\alpha/p} g\big\|_{L^p (\real^n)}
\le 2^{-\alpha/p} \left[  
\frac{\Gamma (\frac{n-\alpha}{2p}) \Gamma (\frac{n}{2p'})}
{\Gamma (\frac{n}{2p'} +\frac{\alpha}{2p}) \Gamma (\frac{n}{2p})}
\right] \big\| (-\Delta)^{\alpha/2p} g\big\|_{L^p (\real^n)}
\end{equation}
\end{cor}

\noindent 
This result recovers the classical Hardy-Rellich inequality for $1<p<n/2$
\begin{equation*}
\int_{\real^n} |x|^{-2p} |g|^p\, dx 
\le \left[\frac{p\, p'}{n(n-2p)}\right]^p
\int_{\real^n} |\Delta g|^p\,dx\ ,
\end{equation*}

\begin{cor}\label{cor2}
{\bf (logarithmic uncertainty)}
For $g\in \S(\real^n)$, $1<p< \infty$
\begin{equation}\label{eq10}
\int_{\real^n}\ln |x| |g|^p\, dx 
+ \int_{\real^n}(\ln\sqrt{-\Delta})g \, g |g|^{p-2}\, dx
\ge D\, \int_{\real^n}  |g|^p\, dx
\end{equation}
\begin{equation*}
D = \ln 2 + \frac{1}{2}\left[ \psi \Big(\frac{n}{2p}\Big) + 
\psi \Big(\frac{n}{2p'}\Big)\right]
\end{equation*}
where $\psi = (\ln\Gamma)'$.
\end{cor}

\noindent
This result follows directly from equation (9) by a differentiation
argument since the constant is identically one with equality at 
$\alpha = 0$. Moreover, it extends the uncertainty principle derived
in \cite4.

The full Stein-Weiss inequality has a more general character which includes 
varying powers of weights (see Appendix). 
Both \cite{17} and \cite{23} suggest that the diagonal map on $L^p (\real^n)$ 
has more difficult aspects, but in fact the inherent convolution structure 
ensures an immediate sharp reduction to radial functions, and the lack of 
extremal functions facilitates the calculation of optimal constants. 

\begin{thm}\label{thm2}
For $f\in L^p(\real^n)$ with $1<p<\infty$, $0<\lambda <n$, $\alpha <n/p$, 
$\beta <n/p'$ and $n= \lambda +\alpha +\beta$ 
\begin{gather}
\Big\| \, |x|^{-\alpha} \Big( |x|^{-\lambda} * \big( |x|^{-\beta} f\big)\Big)
\Big\|_{L^p (\real^n)} 
\le D_{\alpha,\beta}  \|f\|_{L^p (\real^n)}
\label{eq11}\\
\noalign{\vskip6pt}
D_{\alpha,\beta} =    
\pi^{n/2} 
\left[ \frac{\Gamma (\frac{\alpha+\beta}2) 
\Gamma (\frac{n}{2p} - \frac{\alpha}2)
\Gamma (\frac{n}{2p'} - \frac{\beta}2)}
{\Gamma (\frac{n-\alpha-\beta}2) 
\Gamma (\frac{n}{2p'} + \frac{\alpha}2) 
\Gamma (\frac{n}{2p} + \frac{\beta}2)} \right]
\notag
\end{gather}
\end{thm}

\begin{proof} 
Note that it is not required that both $\alpha$ and $\beta$ are non-negative,  
but rather $0<\alpha +\beta <n$. 
By setting $t= |x|$ and 
$\sigma = \frac{n}p - \alpha -\frac{\lambda}2$, 
inequality~\eqref{eq11} has an equivalent formulation as a convolution 
estimate on the product manifold $\real_+ \times S^{n-1}$: 
\begin{equation*}
\| K_{\sigma,\lambda} * h \|_{L^p (\real_+ \times S^{n-1})} 
\le \|K_{\sigma,\lambda} \|_{L^1 (\real_+ \times S^{n-1})} 
\|h\|_{L^p (\real_+ \times S^{n-1})}
\end{equation*}
where $h(t,\xi) = |x|^{n/p} f(|x|,\xi)$ and 
$K_{\sigma,\lambda} (t,\xi\cdot\eta) =t^\sigma (t+\frac1t - 2\xi \cdot\eta)
^{-\lambda/2}$ with $\xi,\eta\in S^{n-1}$.
Then 
\begin{equation*}
\begin{split}
D_{\alpha,\beta} & = \| K_{\sigma,\lambda}\|_{L^1(\real_+,S^{n-1})} 
= \int_{\real_+ \times S^{n-1}} \mkern-38mu
t^{\sigma\,+ \, \lambda/2} (t^2 + 1-2t\xi_1)^{-\lambda/2} \,
d\omega \,\frac{dt}t\\
\noalign{\vskip6pt}
& = \int_{\real^n} |x-\xi|^{-\lambda}  |x|^{-(\alpha \,+\, n/p')}\,dx\\
\noalign{\vskip6pt}
& = \pi^{n/2} \left[ \frac{\Gamma (\frac{\alpha+\beta}2) 
\Gamma (\frac{n}{2p} - \frac{\alpha}2)
\Gamma (\frac{n}{2p'} - \frac{\beta}2)}
{\Gamma (\frac{n-\alpha-\beta}2) 
\Gamma (\frac{n}{2p'} + \frac{\alpha}2) 
\Gamma (\frac{n}{2p} + \frac{\beta}2)} \right]
\end{split}
\end{equation*}
where $\xi$ is a unit vector with $\xi_1$ its first component and 
$d\omega$ denotes surface measure on $S^{n-1}$.

Considerable interest has been given recently to the study of Hardy, Rellich
and weighted Sobolev inequalities on complete Riemannian manifolds with
optimal constants and remainder terms (for example, see \cite1, \cite{10}, 
\cite{13}). Davies and Hinz proved the following two $L^p$ inequalities on
$\real^n$ (see Theorems 12 and 13 in \cite{10} ) which
can be extended using the results and methods for the Stein-Weiss inequality
from Theorem 2 above:
\begin{equation}\label{eq12}
\int_{\real^n} |x|^{-\gamma} |\Delta^m u|^p\, dx 
\ge A_{\gamma, m}\,
\int_{\real^n}|x|^{-\gamma - 2mp} |u|^p\,dx\ ,
\end{equation}

\begin{equation}\label{eq13}
\int_{\real^n} |x|^{-\gamma} |\nabla\Delta^m u|^p\, dx 
\ge B_{\gamma, m}\,
\int_{\real^n}|x|^{-\gamma - (2m+1)p} |u|^p\,dx\ 
\end{equation}
for suitable restrictions on $\gamma$ and $m$.

\begin{thm}\label{thm3}
For $g\in \S(\real^n)$ with $1<p<\infty$, $0<\alpha + \beta <n$, $\alpha <n/p$, 
and $\beta <n/p'$ 
\begin{gather}
\||x|^{-\alpha}g\|_{L^p (\real^n)}
\le 2^{-(\alpha + \beta)} 
\left[ \frac{
\Gamma (\frac{n}{2p} - \frac{\alpha}2)
\Gamma (\frac{n}{2p'} - \frac{\beta}2)}
{\Gamma (\frac{n}{2p'} + \frac{\alpha}2) 
\Gamma (\frac{n}{2p} + \frac{\beta}2)} \right]
\| \, |x|^{\beta} (-\Delta)^{(\alpha + \beta)/2}g
\|_{L^p (\real^n)} 
\label{eq14}
\end{gather}
\end{thm}

\begin{proof}
This result follows directly from Theorem 2. To obtain equation (12), set
$\beta = -\gamma /p$, $\alpha = -\gamma /p + 2m$; then $\alpha + \beta = 2m$
and
$$A_{\gamma,m} = 2^{-2mp}
\left[ \frac{
\Gamma (\frac{n-\gamma}{2p} - m)
\Gamma (\frac{n}{2p'} + \frac{\gamma}{2p})}
{\Gamma (\frac{n}{2p'} + \frac{\gamma}{2p} + m) 
\Gamma (\frac{n-\gamma}{2p} )} \right]^p$$
\renewcommand{\qed}{}
\end{proof}

\begin{thm}\label{thm4}
{\bf (Stein-Weiss potentials with gradient estimates)}
For $f, g\in \S(\real^n)$, $1<p<\infty$, $0<\lambda <n$, $\alpha < n/p$, 
$1-n/p < \beta < n/p'$ and $n+1 = \lambda + \alpha + \beta + \sigma $
\begin{gather}
\big\| \, |x|^{-\alpha} \big( |x|^{-\lambda} * (|x|^{-\sigma}f) \big)
\big\|_{L^p (\real^n)} 
\le F_{\alpha,\beta,\sigma}\,   \big\||x|^{\beta}\nabla f \big\|_{L^p (\real^n)}
\label{eq15}\\
\noalign{\vskip6pt}
F_{\alpha,\beta,\sigma} =    
\pi^{n/2} (\frac{n}p + \beta -1)^{-1} 
\left[ \frac{\Gamma (\frac{\alpha+\beta+\sigma -1}2) 
\Gamma (\frac{n}{2p} - \frac{\alpha}2)
\Gamma (\frac{n}{2p'} - \frac{\beta+\sigma -1}2)}
{\Gamma (\frac{n+1 -\alpha-\beta-\sigma}2) 
\Gamma (\frac{n}{2p'} + \frac{\alpha}2) 
\Gamma (\frac{n}{2p} + \frac{\beta+\sigma -1}2)} \right]
\notag
\end{gather}
\begin{gather}
\big\||x|^{-\alpha}g\big\|_{L^p (\real^n)}
\le \frac{2^{-(\alpha + \beta -1)}}{n/p +\beta -1} 
\left[ \frac{
\Gamma (\frac{n}{2p} - \frac{\alpha}2)
\Gamma (\frac{n}{2p'} - \frac{\beta +\sigma -1}2)}
{\Gamma (\frac{n}{2p'} + \frac{\alpha}2) 
\Gamma (\frac{n}{2p} + \frac{\beta +\sigma -1}2)} \right]
\big\| \, |x|^{\beta} \nabla \big(|x|^{\sigma} (-\Delta)^{(\alpha + \beta -1)/2}g
\big)
\big\|_{L^p (\real^n)} 
\label{eq16}
\end{gather}
\end{thm}

\begin{proof} Note that $1<\alpha + \beta + \sigma < n+1$. The argument can be
reduced to radial functions by setting
$$h(|x|) = \Big[\int_{S^{n-1}} |f(x)|^p \, d\omega \Big]^{1/p} \, , \,
\, |\nabla h| \le \Big[\int_{S^{n-1}} |\nabla f|^p \, d\omega \Big]^{1/p} ;$$
and then applying Young's inequality on $S^{n-1}$ to obtain
$$\big\| \, |x|^{-\alpha} \big( |x|^{-\lambda} * (|x|^{-\sigma}f) \big)
\big\|_{L^p (\real^n)}
\le \big\| \, |x|^{-\alpha} \big( |x|^{-\lambda} * (|x|^{-\sigma}h) \big)
\big\|_{L^p (\real^n)}
$$
$$\le F_{\alpha,\beta,\sigma}\,   \big\||x|^{\beta}\nabla h \big\|_{L^p (\real^n)}
\le F_{\alpha,\beta,\sigma}\,   \big\||x|^{\beta}\nabla f \big\|_{L^p (\real^n)}.
$$
To obtain the inner inequality, set $t=|x|$, $u = |x|^{-\beta -n/p}|h'|$ and observe
that the resulting inequality is now a convolution inequality on the multiplicative
group $\real_{+}$
$$\|K_1 * K_2 * u\|_{L^p(\real_{+})}\, \le \|K_1 * K_2\|_{L^1(\real_{+})}\,
\|u\|_{L^p(\real_{+})}
$$
$$F_{\alpha,\beta,\sigma} = \|K_1 * K_2\|_{L^1(\real_{+})} = \|K_1\|_{L^1(\real_{+})} \,
\|K_2\|_{L^1(\real_{+})}
$$
with
$$K_1(t) =   t^{n/p - \alpha -\lambda/2}
\int_{S^{n-1}} \Big( t + \frac1t - 2\xi_1\Big)^{-\lambda/2} \,d\omega ,
$$
$$K_2(t) = t^{n/p +\beta -1} \, \chi_{[0.1]}(t), \, \|K_2\|_{L^1(\real_{+})} =
(\frac{n}p +\beta -1)^{-1}
$$
$$\|K_1\|_{L^1(\real_{+})} = \int_0^\infty K_1(t)\, \frac1t dt = 
\int_{\real^n} |x-\xi|^{-\lambda}  |x|^{-(\alpha \,+\, n/p')}\,dx
$$
$$ = \pi^{n/2}
\left[ \frac{\Gamma (\frac{\alpha+\beta+\sigma -1}2) 
\Gamma (\frac{n}{2p} - \frac{\alpha}2)
\Gamma (\frac{n}{2p'} - \frac{\beta+\sigma -1}2)}
{\Gamma (\frac{n+1 -\alpha-\beta -\sigma}2) 
\Gamma (\frac{n}{2p'} + \frac{\alpha}2) 
\Gamma (\frac{n}{2p} + \frac{\beta+\sigma -1}2)} \right]
$$
so that
$$F_{\alpha,\beta,\sigma} =    
\pi^{n/2} (\frac{n}p + \beta -1)^{-1} 
\left[ \frac{\Gamma (\frac{\alpha+\beta+\sigma -1}2) 
\Gamma (\frac{n}{2p} - \frac{\alpha}2)
\Gamma (\frac{n}{2p'} - \frac{\beta+\sigma -1}2)}
{\Gamma (\frac{n+1 -\alpha -\beta -\sigma}2) 
\Gamma (\frac{n}{2p'} + \frac{\alpha}2) 
\Gamma (\frac{n}{2p} + \frac{\beta+\sigma -1}2)} \right].$$
To obtain equation (13), set $\alpha = \gamma /p +2m +1$, $\beta = -\gamma /p$
and $\sigma = 0$. Then $\alpha + \beta -1 = 2m$ and
$$B_{\gamma, m} = 2^{-(2m+1)p}
\left[ \frac{
\Gamma (\frac{n-\gamma}{2p} - m-\frac{1}2)
\Gamma (\frac{n}{2p'} + \frac{\gamma}{2p}+\frac{1}2)}
{\Gamma (\frac{n}{2p'} + \frac{\gamma}{2p} + m+\frac{1}2) 
\Gamma (\frac{n-\gamma}{2p}+\frac{1}2 )} \right]^p
$$
\renewcommand{\qed}{}
\end{proof}

Problems with mixed homogeneity have been of interest in relation to embedding
estimates for differential operators. Motivated by earlier work (see \cite5,
\cite6, \cite8) centered on Grushin symmetry and Carnot geometry, and the more 
recent work by Barbatis \cite1, corresponding estimates are obtained for
Stein-Weiss potentials on Euclidean product manifolds with differing homogeneity.
Because of the concentration phenomena that occurs for the sharp $L^1$ Young's 
inequality, one can apply either Euclidean analysis or the underlying hyperbolic 
symmetry to obtain optimal estimates.

\begin{thm}\label{thm5}
{\bf (Stein-Weiss potentials with mixed homogeneity)} For $f\in 
L^p(\real^{n+m})$, $w = (x,v)\in \real^n \times \real^m$, $0 < \lambda < n+m$,
$\alpha < m/p$, $\beta < m/p'$ and $n+m = \lambda + \alpha + \beta$
\begin{gather}
\big\| \, |v|^{-\alpha} \big( |w|^{-\lambda} * (|v|^{-\beta}f) \big)
\big\|_{L^p (\real^{n+m})} 
\le G_{\alpha,\beta}\,   \big\| f \big\|_{L^p (\real^{n+m})}
\label{eq17}\\
\noalign{\vskip6pt}
G_{\alpha,\beta} =    
\pi^{(n+m)/2} 
\left[ \frac{\Gamma (\frac{\alpha+\beta}2) 
\Gamma (\frac{m}{2p} - \frac{\alpha}2)
\Gamma (\frac{m}{2p'} - \frac{\beta}2)}
{\Gamma (\frac{n+m -\alpha-\beta}2) 
\Gamma (\frac{m}{2p'} + \frac{\alpha}2) 
\Gamma (\frac{m}{2p} + \frac{\beta}2)} \right]
\notag
\end{gather}
\end{thm}

\begin{proof}
Apply the sharp $L^1$ Young's inequality for convolution on $\real^n$
followed by the Stein-Weiss argument from Theorem 2 for the $\real^m$
variable. That is,
\begin{equation*}
\begin{split}
\big\| \, |v|^{-\alpha} \big( |w|^{-\lambda} * (|v|^{-\beta}f) \big)
\big\|_{L^p (\real^{n+m})} 
& \le \big\| \, |v|^{-\alpha} \big( J * (|v|^{-\beta}u) \big)
\big\|_{L^p (\real^{m})}\\
\noalign{\vskip6pt}
& = \pi^{n/2} \left[ \frac{\Gamma (\frac{(m - \alpha - \beta)}2)}
{\Gamma (\frac{n+m -\alpha -\beta}2)}\right] \, 
\big\| \, |v|^{-\alpha} \big( |v|^{-\lambda'} * (|v|^{-\beta}u) \big)
\big\|_{L^p (\real^{m})}\\
\noalign{\vskip6pt}
& \le \pi^{(n+m)/2} 
\left[ \frac{\Gamma (\frac{\alpha+\beta}2) 
\Gamma (\frac{m}{2p} - \frac{\alpha}2)
\Gamma (\frac{m}{2p'} - \frac{\beta}2)}
{\Gamma (\frac{n+m -\alpha-\beta}2) 
\Gamma (\frac{m}{2p'} + \frac{\alpha}2) 
\Gamma (\frac{m}{2p} + \frac{\beta}2)} \right] \, \|u\|_{L^p (\real^{m})}\\
\noalign{\vskip6pt}
& = \pi^{(n+m)/2} 
\left[ \frac{\Gamma (\frac{\alpha+\beta}2) 
\Gamma (\frac{m}{2p} - \frac{\alpha}2)
\Gamma (\frac{m}{2p'} - \frac{\beta}2)}
{\Gamma (\frac{n+m -\alpha-\beta}2) 
\Gamma (\frac{m}{2p'} + \frac{\alpha}2) 
\Gamma (\frac{m}{2p} + \frac{\beta}2)} \right] \, \|f\|_{L^p (\real^{n+m})}
\end{split}
\end{equation*}
where $\lambda' = \lambda - n = m - \alpha - \beta$ and
$$u(v) = \Big[\int_{\real^{n}} |f(x,v)|^p \, dx \Big]^{1/p} \, , \,
\, \|u\|_{L^p(\real^m)} \, = \, \|f\|_{L^p(\real^{n+m})}  $$
with
$$J(v) = \int_{\real^n} (x^2 + v^2)^{-\lambda/2}\, dx = |v|^{-(\lambda -n)} \,
\frac{\pi^{n/2}}{\Gamma(n/2)} \, \int_0^\infty t^{n/2 -1} \, (1+t)^{-\lambda/2} 
\, dt.$$
An alternate proof can be given by recognizing that the inequality reduces to
radial functions in the $v$-variable and rewriting the inequality to be on
$\real^{n+1}$ with $y = |v| > 0$, and then transferring the resulting inequality
to $(n+1)$-dimensional hyperbolic space $\mathbb{H}^{n+1}$ with left-invariant Haar
measure $d\nu = y^{-n-1}dxdy$. The result follows by applying the sharp $L^1$
Young's inequality for convolution on a non-unimodular group. This curious aspect
of multiplicity of proofs simply reflects the complexity of the symmetry (see
\cite5, \cite8).
\renewcommand{\qed}{}
\end{proof}

\subsection*{Corollary 1}
For $g\in \S(\real^{n+m})$ with $w = (x,v)\in \real^n \times \real^m$, $1<p<\infty$, 
$0<\alpha + \beta <m$, $\alpha <m/p$, 
and $\beta <m/p'$ 
\begin{gather}
\||v|^{-\alpha}g\|_{L^p (\real^{n+m})}
\le 2^{-(\alpha + \beta)} 
\left[ \frac{
\Gamma (\frac{m}{2p} - \frac{\alpha}2)
\Gamma (\frac{m}{2p'} - \frac{\beta}2)}
{\Gamma (\frac{m}{2p'} + \frac{\alpha}2) 
\Gamma (\frac{m}{2p} + \frac{\beta}2)} \right]
\| \, |v|^{\beta} (-\Delta)^{(\alpha + \beta)/2}g
\|_{L^p (\real^{n+m})} 
\label{eq18}
\end{gather}

\noindent Similar arguments to those used for Theorem 5 give gradient estimates with
mixed homogeneity:

\subsection*{Corollary 2}
For $g\in \S(\real^{n+m})$ with $w = (x,v)\in \real^n \times \real^m$, $1<p<\infty$, 
$0<\alpha + \beta <m$, $\alpha <m/p$, 
and $\beta <m/p'$ 
\begin{gather}
\||v|^{-\alpha}g\|_{L^p (\real^{n+m})}
\le 2^{-(\alpha + \beta)} 
\left[ \frac{
\Gamma (\frac{m}{2p} - \frac{\alpha}2)
\Gamma (\frac{m}{2p'} - \frac{\beta}2 +\frac{1}2)}
{\Gamma (\frac{m}{2p'} + \frac{\alpha}2) 
\Gamma (\frac{m}{2p} + \frac{\beta}2 +\frac{1}2)} \right]
\| \, |v|^{\beta} \nabla (-\Delta)^{(\alpha + \beta -1)/2}g
\|_{L^p (\real^{n+m})} 
\label{eq19}
\end{gather}

Riesz potentials can be viewed as extending Young's inequality to include
weak Lebesgue classes.
Stein-Weiss potentials extend Riesz potentials by incorporating successive 
multiplication by fractional powers --- first on the function side and then 
followed on the Fourier transform side with repetitions. 
Such iterations include both extensions of Young's inequality and the 
framework of the Maz'ya-Eilertsen inequality \cite{11}. 
The following theorem extends from $L^2 (\real^n)$ to $L^p(\real^n)$ 
results on iterated Stein-Weiss potentials from \cite7.
\renewcommand{\qed}{}
\end{proof}

\begin{thm}\label{thm6}
For $f\in L^p (\real^n)$, $1<p<\infty$, $\alpha +\sigma = \beta +\rho$,
$0<\alpha < n/p'$, $0<\sigma <n$, $-n/p < \rho-\sigma$, $\rho<n/p'$ 
\begin{gather}
\Big\|\, |x|^{-(n-\alpha)} * \Big[ |x|^{-\beta} (|x|^{-(n-\sigma)} 
* (|x|^{-\rho} f))\Big] \Big\|_{L^p (\real^n)} 
\le E_{\alpha,\sigma,\rho}\|f\|_{L^p(\real^n)}
\label{eq20}\\
\noalign{\vskip6pt}
E_{\alpha,\sigma,\rho} =  \pi^n \left[ 
\frac{\Gamma (\frac{\alpha}2)\, \Gamma (\frac{\sigma}2)\, 
\Gamma (\frac{n}{2p})\, \Gamma (\frac{n}{2p'} -\frac{\alpha}2)\, 
\Gamma (\frac{n}{2p'} - \frac{\rho}2)\, 
\Gamma (\frac{n}{2p} + \frac{\rho}2 - \frac{\sigma}2)} 
{\Gamma (\frac{n-\alpha}2)\, \Gamma (\frac{n-\sigma}2)\, 
\Gamma (\frac{n}{2p'})\, \Gamma (\frac{n}{2p} + \frac{\alpha}2)\, 
\Gamma (\frac{n}{2p} + \frac{\rho}2)\, 
\Gamma (\frac{n}{2p'} - \frac{\rho}2 +\frac{\sigma}2)}
\right]\ .\notag
\end{gather}
\end{thm}

\begin{proof}
Set $\lambda = \frac{n}2 - \frac{n}{p'} + \rho -\frac{\sigma}2$ 
and $\nu = \frac{n}p - \frac{n}2 + \frac{\alpha}2$; 
then inequality~\eqref{eq20} has an equivalent formulation as a 
convolution estimate on the product manifold $\real_+ \times S^{n-1}$: 
\begin{gather*}
\|K_1 * K_2 * h\|_{L^p (\real_+\times S^{n-1})} 
\le \|K_1 * K_2\|_{L^1 (\real_+\times S^{n-1})}
\|h\|_{L^p(\real_+ \times S^{n-1})} \\
\noalign{\vskip6pt}
E_{\alpha,\sigma,\rho} = 
\| K_1 * K_2\|_{L^1(\real_+ \times S^{n-1})} 
= \|K_1\|_{L^1 (\real_+ \times S^{n-1})} 
\|K_2\|_{L^1 (\real_+ \times S^{n-1})}
\end{gather*}
where $h(t,\xi) = |x|^{n/p} f(|x|,\xi)$ with $t= |x|$ and 
\begin{gather*}
K_1(t,\xi\cdot\eta) = t^\nu\Big(t+\frac1t - 2\xi\cdot\eta\Big)^{-(n-\alpha)/2}\\
\noalign{\vskip6pt}
K_2(t,\xi\cdot\eta) = t^\lambda\Big(t+\frac1t -2\xi\cdot\eta\Big)^
{-(n-\sigma)/2}\ .
\end{gather*}
\begin{equation*}
\begin{split}
\|K_1 \|_{L^1(\real_+\times S^{n-1})} 
& = \int_{\real_+\times S^{n-1}} \mkern-38mu
t^{\nu + (n-\alpha)/2} (t^2 + 1 -2t\xi_1)^{-(n-\alpha)/2}\,d\omega\frac{dt}t\\
\noalign{\vskip6pt}
& = \int_{\real^n} |x-\xi|^{-(n-\alpha)} |x|^{-(n+\alpha)/2\, +\,\nu}\,dx\\
\noalign{\vskip6pt}
& = \int_{\real^n} |x-\xi|^{-(n-\alpha)} |x|^{-n/p'}\,dx\\
\noalign{\vskip6pt}
& = \pi^{n/2} \left[ \frac{\Gamma (\frac{\alpha}2)\, 
\Gamma (\frac{n}{2p})\, \Gamma (\frac{n}{2p'} - \frac{\alpha}2)}
{\Gamma (\frac{n-\alpha}2)\, \Gamma (\frac{n}{2p'})\, 
\Gamma (\frac{n}{2p} + \frac{\alpha}2)} 
\right]
\end{split}
\end{equation*}
\begin{equation*}
\begin{split}
\|K_2 \|_{L^1(\real_+\times S^{n-1})} 
& = \int_{\real_+\times S^{n-1}} \mkern-38mu
t^{\lambda +(n-\sigma)/2}(t^2 +1-2t\xi_1)^{-(n-\sigma)/2}\,d\omega\frac{dt}t\\
\noalign{\vskip6pt}
& = \int_{\real^n} |x-\xi|^{-(n-\sigma)} |x|^{-(n+\sigma)/2\, +\,\lambda}\,dx\\
\noalign{\vskip6pt}
& = \int_{\real^n} |x-\xi|^{-(n-\sigma)} |x|^{-n/p' + \rho -\sigma}\, dx\\
\noalign{\vskip6pt}
& = \pi^{n/2} \left[ \frac{\Gamma (\frac{\sigma}2)\, 
\Gamma (\frac{n}{2p} + \frac{\rho}2 - \frac{\sigma}2)\, 
\Gamma (\frac{n}{2p'} - \frac{\rho}2)}
{\Gamma (\frac{n-\sigma}2)\, 
\Gamma (\frac{n}{2p'} +\frac{\sigma}2 -\frac{\rho}2)\, 
\Gamma (\frac{n}{2p} + \frac{\rho}2)} 
\right]\ .
\end{split}
\end{equation*}
Then 
\begin{equation*}
E_{\alpha,\sigma,\rho} =  \pi^n \left[
\frac{\Gamma (\frac{\alpha}2)\, \Gamma (\frac{\sigma}2)\,
\Gamma (\frac{n}{2p})\, \Gamma (\frac{n}{2p'} -\frac{\alpha}2)\,
\Gamma (\frac{n}{2p'} - \frac{\rho}2)\,
\Gamma (\frac{n}{2p} + \frac{\rho}2 - \frac{\sigma}2)}
{\Gamma (\frac{n-\alpha}2)\, \Gamma (\frac{n-\sigma}2)\,
\Gamma (\frac{n}{2p'})\, \Gamma (\frac{n}{2p} + \frac{\alpha}2)\,
\Gamma (\frac{n}{2p} + \frac{\rho}2)\,
\Gamma (\frac{n}{2p'} - \frac{\rho}2 +\frac{\sigma}2)}
\right]\ .
\end{equation*}
Note that as in the original Stein-Weiss inequality, $\beta$ and $\rho$  
can take negative values as long as $\beta +\rho$ is positive. 
Two lemmas are used in the above calculation of the optimal constant. The first
lemma is obtained using the inverse Fourier transform for fractional powers.
\renewcommand{\qed}{}
\end{proof}

\begin{lem}\label{lem1}
For $0<\lambda <n$, $0<\mu <n$ and $\lambda +\mu >n$
\begin{equation}\label{eq21}
|x|^{-\lambda} * |x|^{-\mu} 
= \pi^{n/2} \left[ 
\frac{ \Gamma (\frac{n-\lambda}2)\, \Gamma (\frac{n-\mu}2)\, 
\Gamma (\frac{\lambda+\mu -n}2)}
{\Gamma (\frac{\lambda}2)\, \Gamma (\frac{\mu}2)\, 
\Gamma (\frac{2n -\lambda -\mu}2)} \right]\ \, |x|^{-(\lambda + \mu -n)}.
\end{equation}
\end{lem}

\begin{lem}\label{lem2}
For $F\in L^p (S^n)$, $G\in L^{p'} (S^n)$ and $K(\xi\cdot\eta) \ge0$ 
with $\xi,\eta \in S^n$ 
\begin{gather}
\Big| \int_{S^n\times S^n} G(\xi) K (\xi\cdot\eta) F(\eta)\,d\xi\,d\eta\Big|
\le \bigg(\int_{S^n} K(\xi_1) \,d\xi\bigg) 
\|F\|_{L^p (S^n)} \|G\|_{L^{p'} (S^n)}
\label{eq22}\\
\noalign{\vskip6pt}
\|K * F\|_{L^p(S^n)}  \le \|K\|_{L^1 (S^n)} \|F\|_{L^p (S^n)}
\label{eq23}
\end{gather}
\end{lem}

\begin{proof}
Split the integrand in \eqref{eq22} into the product of two terms, 
$FK^{1/p}$ and $GK^{1/p'}$ and apply H\"older's inequality on the product 
manifold $S^n\times S^n$.
\renewcommand{\qed}{}
\end{proof}

\begin{Cor-Hardy}[iterated Hardy-Rellich inequality]
For $h\in \S (\real^n)$, $1<p<n$, and the parameters $\alpha$, $\beta$, 
$\sigma$, $\rho$ as specified in Theorem~\ref{thm6} 
\begin{gather}
\|h\|_{L^p(\real^n)} 
\le  C\big\|\, |x|^\rho (-\Delta)^{\sigma/2} |x|^\beta (-\Delta)^{\alpha/2}
h \Big\|_{L^p (\real^n)}
\label{eq24}\\
\noalign{\vskip6pt}
C = 2^{- (\alpha +\sigma)} 
\left[ \frac{ \Gamma (\frac{n}{2p})\, 
\Gamma (\frac{n}{2p'} -\frac{\alpha}2)\, 
\Gamma (\frac{n}{2p'} - \frac{\rho}2)\,
\Gamma (\frac{n}{2p} + \frac{\rho}2 - \frac{\sigma}2) }
{\Gamma (\frac{n}{2p'})\, 
\Gamma (\frac{n}{2p} + \frac{\alpha}2)\, 
\Gamma (\frac{n}{2p} + \frac{\rho}2)\, 
\Gamma (\frac{n}{2p'} - \frac{\rho}2 + \frac{\sigma}2)} 
\right]\ .
\notag
\end{gather}
\end{Cor-Hardy}

As with Riesz potentials, iterated Stein-Weiss potentials coupled with 
symmetrization provide extensions of Young's inequality to include 
weak Lebesgue classes: 
\begin{equation*}
h\in L_{q,\infty} (\real^n) \Longleftrightarrow h^* (x) 
\le \frac{A}{|x|^{n/q}}\ ,\qquad 0<q<\infty
\end{equation*}
where $h^*$ denotes the equimeasurable radial decreasing rearrangement of 
$|h|$ on $\real^n$.

\begin{Cor-Young}[Young's inequality]
For $f\in L^p (\real^n)$, 
$1<p<\infty$, $1<q_k<\infty$, $\sum 1/q_k = 2$, $q_1<p$, $q_4 >p'$ 
\begin{equation}\label{eq25}
\| h_1 *h_2 (h_3 *(h_4 f))\|_{L^p(\real^n)}
\le C\prod_{k=1}^4 \|h_k\|_{L_{q_{k,\infty}}(\real^n)} \|f\|_{L^p (\real^n)}
\end{equation}
\end{Cor-Young}

By using the methods outlined in the Appendix (section~3), one obtains 
results from the Stein-Weiss Theorem for Young's inequality that include 
off-diagonal $L^p-L^q$ maps:

\begin{thm}\label{thm7}
For $f\in L^p(\real^n)$, $1<p\le q<\infty$, $1<r_k <\infty$, 
$1/q\, +\, 1/p' =\sum 1/r_k$, $r_1 >q$, $1< q_k <\infty$, 
$1\,+\, 1/q\,+\, 1/p' = \sum 1/q_k$, $q_1 < q$, $q_4 >p'$ 
\begin{gather}
\|g_1 (g_2 *(g_3f))\|_{L^q (\real^n)} 
\le C\prod_{k=1}^3 \|g_k\|_{L_{r_{k,\infty}}(\real^n)} \|f\|_{L^p(\real^n)}
\label{eq26}\\
\noalign{\vskip6pt}
\|h_1 *h_2(h_3 * (h_4 f))\|_{L^q (\real^n)}
\le C \prod_{k=1}^4 \|h_k\|_{L_{q_{k,\infty}} (\real^n)}\|f\|_{Lp(\real^n)}\ .
\label{eq27}
\end{gather}
\end{thm}

\begin{remark}
The constants obtained in Theorems~\ref{thm1}, \ref{thm2}, \ref{thm3}, 
\ref{thm4}, \ref{thm5} and \ref{thm6} 
are optimal and no extremal functions exist for the corresponding variational 
problems. 
This circumstance follows directly from recognition that dilation invariance 
implies an equivalent formulation as sharp convolution estimates on 
$\real_+ \times S^{n-1}$ with a positive integrable kernel. After this paper 
was finished, the author learned that Theorems 2 and 3 had been obtained 
independently by Samko \cite{19}.
\end{remark}

\section*{Appendix}

Three arguments are developed here that contribute to a broader setting
for the results in this note. A short proof of the classical Pitt's
inequality on $\real^n$ for $n\ge 2$ is given by applying rearrangement
inequalities, convolution estimates and interpolation.
Pitt's inequality can be viewed as a generalized Hausdorff-Young
inequality with weights. From that context, it is not surprising that
the only cases where optimal constants have been calculated correspond
to the Hausdorff-Young inequality, the Hardy-Littlewood-Sobolev inequality,
and the spectral-level diagonal maps defined by Stein-Weiss potentials.

The second issue addressed here is to show the connection between an
argument taken from Stein and Weiss [23] and the Hardy-Littlewood paradigm
relating dilation invariance for a positive linear operator to sharp
convolution estimates on the multiplicative group. To fully complete the
background for these problems, a short proof of the Stein-Weiss theorem on
fractional integrals is given using Young's inequality. The structural
correspondence between Pitt's inequality and the Stein-Weiss theorem
underlines the essential duality between the Fourier transform and
convolution.
More broadly, the essential idea is that identification of global symmetry 
facilitates not simply calculation of the mathematical estimates but also 
a fuller understanding of the analytic framework.

\section*{1. Pitt's inequality}

\begin{Pittsinequality}
For $f\in \S(\real^n)$, $1<p\le q<\infty$, $0<\alpha < n/q$, $0<\beta< 
n/p'$
and $n\ge 2$
\begin{equation}\label{eq28}
\bigg[ \int_{\real^n} \big|\, |x|^{-\alpha} \hat f\big|^q\,dx\bigg]^{1/q}
\le A\bigg[ \int_{\real^n} \big|\, |x|^\beta f\big|^p\,dx\bigg]^{1/p}
\end{equation}
with the index constraint
\begin{equation*}
\frac{n}p +\frac{n}q + \beta -\alpha = n
\end{equation*}
(note that primes denote dual exponents, $1/p \,+\, 1/p' =1$).
\end{Pittsinequality}

\begin{proof}
The argument will be accomplished in stages.
\renewcommand{\qed}{}
\end{proof}

\step{1}    
Observe that by duality, inequality \eqref{eq28}  for $p<q$ implies the
dual inequality
\begin{equation}\label{eq29}
\bigg[ \int_{\real^n} \big|\, |x|^{-\beta} \hat h\big|^{p'}\,dx\bigg]^{1/p'}
\le A\bigg[ \int_{\real^n}\big|\, |x|^\alpha h\big|^{q'}\, dx\bigg]^{1/q'}
\end{equation}
which simply reverses the place of $\alpha$ and $\beta$ with
$q' < p'$.

\step{2}    
As a preliminary calculation, consider \eqref{eq28}  for $q=p'$ which
requires $\beta-\alpha =0$.
Then using the Hausdorff-Young inequality combined with Theorem~\ref{thm1}
\begin{equation*}
\begin{split}
\bigg[ \int_{\real^n} \big|\, |x|^{-\alpha}\hat f\big|^q\,dx\bigg]^{1/q}
& \le C_1 \bigg[ \int_{\real^n} \big|\, |x|^{-(n-\alpha)} *f\big|^p\, dx
\bigg]^{1/p}\\
\noalign{\vskip6pt}
& \le C_2 \bigg[ \int_{\real^n} \big|\, |x|^\alpha f\big|^p\,dx\bigg]^{1/p}
\end{split}
\end{equation*}
which demonstrates \eqref{eq28} for the case $q=p'$ and $0<\alpha < n/q$.

\step{3}
For $q=2$ and $\beta=0$, then \eqref{eq28} corresponds to the
Hardy-Littlewood-Sobolev inequality.
But it is surprising to note that the case of positive $\beta$ does not
follow from this result.
For $\beta>0$ and $q=2$, Pitt's inequality is equivalent to
\begin{equation*}
\bigg[ \int_{\real^n} \big|\, |x|^{-\alpha} * (|x|^{-\beta} f)\big|^2\,dx
\bigg]^{1/2}
\le B \bigg[ \int_{\real^n} |f|^p\, dx\bigg]^{1/p}
\end{equation*}
with $\frac{n}2 - \alpha = \frac{n}{p'} -\beta$.
Using symmetrization and rearrangement, this estimate is reduced to
non-negative radial decreasing functions on $\real^n$.
Set $t= |x|$ and $h(t) = |x|^{-n/p} f(x)$.
Now the above inequality is equivalent to
$$\|\psi_\alpha * h\|_{L^2(\real_+)}
\le C\|h\|_{L^p (\real_+)}$$
where
$$\psi_\alpha (t) = t^{\alpha/2} \int_{S^{n-1}} (t+t^{-1}
- 2\xi_1)^{-(n-\alpha)/2}\,d\xi$$
with $d\xi$ being normalized surface measure on $S^{n-1}$, $\xi_1$
denoting the first component of $\xi$ and Haar measure on $\real_+$
being $t^{-1}\,dt$.
Since $\psi_\alpha \in L^r (\real_+)$ for all $r\ge1$, inequality 
\eqref{eq28}
for this case follows from Young's inequality
$$\| \psi_\alpha * h \|_{L^2 (\real_+)}
\le \|\psi_\alpha \|_{L^r(\real_+)} \|h\|_{L^p(\real_+)}$$
for $1/2 = 1/p\, +\, 1/r \, -\, 1$.

\step{4}
The central step in the proof is to show that \eqref{eq28} holds whenever
$q$ is an even integer, $q=2m$ with $m\ge2$.
Then
\begin{gather*}
\int_{\real^n} \big|\, |x|^{-\alpha} \hat f \big|^{2m}\, dx
= c \int_{\real^n} \big| \big(|x|^{-(n-\alpha)} *f\big) *\cdots *
\big( |x|^{-(n-\alpha)} * f\big)\big|^2\,dx\\
\noalign{\vskip6pt}
\le c\int_{\real^n} \big| \big( |x|^{-(n-\alpha)} * f^*\big) * \cdots *
\big( |x|^{-(n-\alpha)} * f^*\big)\big|^2\,dx
\end{gather*}   
where $f^*$ denotes the radial equimeasurable decreasing
rearrangement of $|f|$, and there are $m$ convolutions in each
string appearing in the integral above,  
and for any $p>0$ and  $\alpha >0$
\begin{equation*}
\int_{\real^n} \big|\, |x|^\alpha f^*\big|^p \,dx
\le \int_{\real^n} \big|\, |x|^\alpha f\big|^p\,dx
\end{equation*}
since $|x|^\alpha$ is an increasing function.
Hence for $q=2m$, it suffices to consider \eqref{eq28} for non-negative
radial decreasing functions.
The most direct approach at this point is to apply Young's inequality
successively --- first on $\real^n$ and then on $\real_+$.
Inequality \eqref{eq28} is equivalent to
\begin{equation*}
\bigg[ \int_{\real^n} \big|\, |x|^{-\alpha} \big( |x|^{-\beta} 
f\big)^\wedge\, \big|^{2m}
\,dx\bigg]^{1/2m}
\le A \bigg[ \int_{\real^n} |f|^p\,dx \bigg]^{1/p}
\end{equation*}
and
\begin{equation*}
\begin{split}
\bigg[\int_{\real^n} \big|\, |x|^{-\alpha} \big( |x|^{-\beta} 
f\big)^\wedge\,
\big|^{2m}\,dx\bigg]^{1/2m}
&= c \bigg[ \int_{\real^n} \big| \big( |x|^{-(n-\alpha)} *  |x|^{-\beta} 
f\big)
* \cdots * \big( |x|^{-(n-\alpha)} * |x|^{-\beta} 
f\big)\big|^2\,dx\bigg]^{1/2m}\\
\noalign{\vskip6pt}
&\le c \big\|\, |x|^{-(n-\alpha)} * |x|^{-\beta} f\big\|_{L^r(\real^n)}
\end{split}
\end{equation*}
by applying Young's inequality on $\real^n$ for $r= 2m/(2m-1)$.
The proof for this case is finished by showing that
\begin{equation*}
\big\| \, |x|^{-(n-\alpha)} * |x|^{-\beta} f\big\|_{L^r (\real^n)}
\le C_1 \| f\|_{L^p (\real^n)}
\end{equation*}
with $n/r = (\beta-\alpha) + n/p$.
It suffices to show this inequality for radial functions $f$, and that case
is controlled by Young's inequality on $\real_+$
\begin{equation*}
\|\psi_{\sigma,\alpha} * h\|_{L^r(\real_+)}
\le \|\psi_{\sigma,,\alpha}\|_{L^s(\real_+)} \|h\|_{L^p(\real_+)}
\end{equation*}
where $\frac1r = \frac1s +  \frac1p - 1$, $r=2m/(2m-1)$,
$\sigma = \frac{n}2 - \frac{n}{2m} +\frac{\alpha}2$ and
\begin{equation*}
\psi_{\sigma,\alpha} = t^\sigma \int_{S^{n-1}} (t+t^{-1} - 2\xi_1)
^{-(n-\alpha)/2}\,d\xi
\end{equation*}
with $d\xi$ being normalized surface measure on $S^{n-1}$ and
$\xi_1$ the first component of $\xi$.
$\psi_{\sigma,\alpha} \in L^s (\real_+)$ for any $s\ge 1$ if
$(n-\alpha)/2 -\sigma = n/2m \, -\, \alpha \, =\, n/q\,-\, \alpha >0$
which is a necessary restriction for Pitt's inequality.
Another approach for $m\ge 2$ would be to use the fact that
the Fourier transform of a radial function results in a transform
with a kernel having some decay at infinity:
\begin{equation*}
\int_{\real^n} e^{2\pi ixy} f(|y|)\, dy
= c \int_0^\infty (|x|\, |y|)^{-(n/2\,-\,1)} J_{n/2\,-\,1}
(2\pi |x|\, |y|) f(|y|)\,dy \ .
\end{equation*}
Here $J_{n/2\,-\,1}$ is a Bessel function of the first kind.
Then inequality \eqref{eq28} for radial functions
\begin{equation*}
\big\|\, |x|^{-\alpha} \big( |x|^{-\beta} f\big)^\wedge
\big\|_{L^q(\real^n)} \le A \|f\|_{L^p (\real^n)}
\end{equation*}
corresponds to the convolution inequality on $\real_+$
\begin{equation*}
\| K*h\|_{L^q(\real_+)} \le B \|h\|_{L^p (\real_+)}
\end{equation*}
with $K(t) = t^{\sigma + 1\, - \, n/2} J_{n/2\,-\,1} (t)$
$\sigma = \frac{n}q - \alpha =\frac{n}{p'} -\beta >0$.
If $q$ is at least 4, then $K\in L^r (\real_+)$ for all $r\ge1$.
By Young's inequality, one can take $B = \|K\|_{L^r(\real_+)}$
with $\frac1r =  \frac1q +\frac1{p'}$.

\step{5}
The proof of Pitt's inequality is completed by applying
Riesz-Stein complex interpolation to show that the linear operator
\begin{equation*}
(T_{\alpha,\beta} f) (x) = |x|^{-\alpha} 
\big(|x|^{-\beta}f\big)^\wedge\, (x)
\end{equation*}
is bounded from $L^p (\real^n)$ to $L^q(\real^n)$ on the convex
domain $\frac{n}q - \alpha= \frac{n}{p'} -\beta$,
$q\ge p$, $q\ge 2$, $0<\alpha < \frac{n}q$, and $0<\beta <\frac{n}{p'}$.
The operator $T_{\alpha,\beta}$ is an analytic function of both $\alpha$
and $\beta$.
Every point in the above set lies on a line connecting two values
of $q$ that are even integers.
Using Steps~3 and 4 which give the required end-point bounds,
then Pitt's inequality holds for all parameters in the above
3-dimensional convex domain (see lemma on page 385 in \cite{22}.
Duality from Step~1 extends the result to values $p<q<2$ and the
proof of \eqref{eq28} is finished.

\section*{2. Stein-Weiss lemma}

The Hardy-Littlewood paradigm relating dilation invariance for
a positive integral operator to sharp convolution estimates on the
multiplicative group is used to clarify an old argument by
Stein and Weiss (see Lemma~2.1 in \cite{23}).

\begin{SW}    
Suppose $K$ is a non-negative kernel defined on $\real^n\times
\real^n$, continuous on any domain that excludes the point $(0,0)$,
homogeneous of degree $-n$, $K(\delta u,\delta v) = \delta^{-n} K(u,v)$,
and $K(Ru,Rv) = K(u,v)$ for any $R\in SO(n)$.
Then $K$ defines an integral operator
\begin{equation*}
(Tf) (x) = \int_{\real^n} K(x,y) f(y)\,dy
\notag
\end{equation*}
which maps $L^p (\real^n)$ to $L^p(\real^n)$ for $1<p<\infty$
\begin{equation}\label{eq30}
\|Tf\|_{L^p(\real^n)} \le A\|f\|_{L^p(\real^n)}
\end{equation}
where the optimal constant is given by
\begin{equation*}
A = \int_{\real^n} K(x,\hat e_1) |x|^{-n/p'}\,dx
\end{equation*}
and $\hat e_1$ is a unit vector in the first coordinate direction.
\end{SW}

\begin{proof}
Set $s = |y|$, $t= |x|$, $h(s) = s^{n/p} [\int_{S^{n-1}} 
|f(s\xi)|^p\,d\xi]^{1/p}$
and $\bar K (t,s) = \int_{S^{n-1}} K(t\xi,s\hat e_1)\,d\sigma$.
Then using the invariance of $K$ with respect to the rotation group
\begin{equation*}
\begin{split}
\|Tf\|_{L^p(\real^n)}
&   \le \bigg[ \int_0^\infty \Big| t^{n/p} \int_0^\infty \bar K
(t,s) s^{n/p'} h(s) \frac{ds}{s} \Big|^p \frac{dt}t\bigg]^{1/p}\\
\noalign{\vskip6pt}
&= \bigg[ \int_0^\infty \Big| \int_0^\infty \Big(\frac{t}s\Big)^{n/p}
\bar K\Big(\frac{t}s,1\Big) h(s) \frac{ds}s \Big|^p \frac{dt}t 
\bigg]^{1/p}\\
\noalign{\vskip6pt}
& = \|\psi * h\|_{L^p (\real_+)}
\le \|\psi\|_{L^1 (\real_+)} \|h\|_{L^p(\real_+)}
\end{split}
\end{equation*}
where $\psi (t) = t^{n/p} \bar K(t,1)$,
$\|h\|_{L^p(\real_+)}= \|f\|_{L^p(\real^n)}$ and
\begin{equation*}
\|\psi\|_{L^1(\real_+)} = \int_{\real^n} K(x,\hat e_1)
|x|^{-n/p'}\,dx\ .
\end{equation*}
One can see that $\|\psi\|_{L^1 (\real_+)}$ is the optimal constant
by considering radial functions.
\renewcommand{\qed}{}
\end{proof}

\section*{3. Stein-Weiss theorem}

\begin{SWthm} 
For $f\in \S(\real^n)$, $1<p\le q<\infty$, $0<\lambda <n$, $\alpha < n/q$, 
$\beta < n/p'$, $\alpha +\beta \ge 0$ and $\frac{n}q + \frac{n}{p'} = 
\lambda +\alpha +\beta$
\begin{equation}\label{eq31} 
\big\|\, |x|^{-\alpha} \big( |x|^{-\lambda} * f\big) \big\|_{L^q(\real^n)} 
\le A\big\| \, |x|^\beta f\big\|_{L^p (\real^n)}\ .
\end{equation}
\end{SWthm}

\begin{proof}
Observe that because $\alpha$ and $\beta$ are not both necessarily 
non-negative but only the sum is non-negative, there is not a reduction 
to radial functions. 
But it is still possible to express this problem as a product convolution 
on $\real_+ \times S^{n-1}$ where convolution on the symmetric space $S^{n-1}$
is defined by 
$$(K*F) (\xi) = \int_{S^{n-1}} K(\xi\cdot n) F(\eta)\,d\eta$$
and $K$ is defined on the interval $(-1,1)$. 
The standard convolution estimates expressed by W.H.~Young's inequality 
hold in this symmetric space setting (see Zygmund \cite{25}: Section~1.15, ch.2, 
vol.I). 
Then inequality~\eqref{eq31} is determined by 
\begin{equation}\label{eq32} 
\|K_\sigma * h\|_{L^q (\real_+\times S^{n-1})} 
\le \|K_\sigma\|_{L^r (\real_+ \times S^{n-1})} 
\|h\|_{L^p(\real_+\times S^{n-1})}
\end{equation}
where $h(t,\xi) = |x|^{-\beta -\, n/p} f(|x|,\xi)$ with $t = |x|$, 
$\sigma = \frac{n}q - \alpha - \frac{\lambda}2 = - (\frac{n}{p'} -\beta 
- \frac{\lambda}2)$, and $1/r\,=\, 1/q\,+\, 1/p' = (\lambda +\alpha 
+\beta)/n$. 
To obtain the Stein-Weiss theorem, one need only check that 
$$K_\sigma (t,\xi\cdot\eta) = t^\sigma \Big(t+\frac1{t} - 2\xi \cdot\eta
\Big)^{-\lambda/2} \in L^r (\real_+ \times S^{n-1})$$
where $r = n/(\lambda + \alpha +\beta)$. 
Note that $r\lambda <n$.
Calculate 
$$\int_{\real_+ \times S^{n-1}} 
|K_\sigma (t,\xi_1)|^r d\omega \frac{dt}{t} 
= \int_{\real^n} \frac1{|x-\eta|^\lambda} \ \frac1{|x|^{\alpha+n/q'}} dx 
<\infty$$
if $0<\lambda<n$ and $-n/q'< \alpha < n/q$ which follows from 
$-n/q' < -n/p' <-\beta <\alpha <n/q$, and this completes the proof
of the Stein-Weiss theorem.
\renewcommand{\qed}{}
\end{proof}

\section*{Acknowledgement}

I would like to thank Michael Perelmuter for drawing my attention to the
work of Herbst \cite{15} and the paper \cite{16} where the spectral properties
of Schr\"odinger operators acting on $L^p$ spaces are studied.


\end{document}